\documentclass[10pt,a4paper]{article}
\usepackage[margin=1.1cm]{geometry}
\usepackage{amsmath,amsthm,amsfonts,amssymb,amscd,cite,graphicx}

\usepackage{titlesec}
\titleformat{\section}
{\normalfont\fontsize{12}{15}\bfseries}{\thesection}{1em.}{}

\allowdisplaybreaks[4]

\let\oldbibliography\thebibliography
\renewcommand{\thebibliography}[1]{%
  \oldbibliography{#1}%
  \setlength{\itemsep}{-2pt}%
}

\begin{document}

\baselineskip=0.20in

\newcommand{\tbox}[1]{\mbox{\tiny #1}}

\makebox[\textwidth]{%
\hglue-15pt
\begin{minipage}{0.6cm}	
\vskip9pt
\end{minipage} \vspace{-\parskip}
\begin{minipage}[t]{6cm}
\footnotesize{ {\bf Discrete Mathematics Letters} \\ \underline{www.dmlett.com}}
\end{minipage}
\hfill
\begin{minipage}[t]{6.5cm}
\normalsize {\it Discrete Math. Lett.}  {\bf X} (202X) XX--XX
\end{minipage}}
\vskip36pt

\noindent
{\large \bf Stolarsky--Puebla index}\\

\noindent
J. A. M\'endez-Berm\'udez$^{1,}\footnote{Corresponding author (jmendezb@ifuap.buap.mx)}$, R. Aguilar-S\'anchez$^{2}$,  Ricardo Abreu Blaya$^{3}$, Jos\'e M. Sigarreta$^{3}$\\

\noindent
\footnotesize 
$^1${\it Instituto de F\'{\i}sica, Benem\'erita Universidad Aut\'onoma de Puebla, Apartado Postal J-48, Puebla 72570, Mexico} \\
\noindent
$^2${\it Facultad de Ciencias Qu\'imicas, Benem\'erita Universidad Aut\'onoma de Puebla,
Puebla 72570, Mexico}\\
\noindent
$^3${\it Facultad de Matem\'aticas, Universidad Aut\'onoma de Guerrero, Carlos E. Adame No.54 Col. Garita, Acapulco Gro. 39650, Mexico} \\

\noindent
 (\footnotesize Received: Day Month 202X. Received in revised form: Day Month 202X. Accepted: Day Month 202X. Published online: Day Month 202X.)\\

\setcounter{page}{1} \thispagestyle{empty}

\baselineskip=0.20in

\normalsize

 \begin{abstract}
 \noindent
 We introduce a degree--based variable topological index inspired on the Stolarsky mean (known 
as the generalization of the logarithmic mean). 
We name this new index as the Stolarsky--Puebla index: 
$SP_\alpha(G) = \sum_{uv \in E(G)} d_u$, if $d_u=d_v$, and
$SP_\alpha(G) = \sum_{uv \in E(G)} \left[\left( d_u^\alpha-d_v^\alpha\right)/\left( \alpha(d_u-d_v\right)\right]^{1/(\alpha-1)}$, otherwise.
Here, $uv$ denotes the edge of the network $G$ connecting the vertices $u$ and $v$,
$d_u$ is the degree of the vertex $u$, and $\alpha \in \mathbb{R} \backslash \{0,1\}$.
Indeed, for given values of $\alpha$, the Stolarsky--Puebla index reproduces well-known topological indices
such as the reciprocal Randic index, the first Zagreb index, and several mean Sombor indices.
Moreover, we apply these indices to random networks and demonstrate that 
$\left< SP_\alpha(G) \right>$, normalized to the order of the network, scale with 
the corresponding average degree $\left< d \right>$.
 \\[2mm]
 \end{abstract}

\baselineskip=0.20in

\section{Introduction}

For two positive real numbers $x$, $y$ the Stolarsky mean $S_\alpha(x,y)$ is defined as~\cite{S75}
\begin{equation}
S_\alpha(x,y)=\lim_{(\xi,\eta)\to(x,y)} \left( \frac{\xi^\alpha-\eta^\alpha}{\alpha(\xi-\eta)} \right)^{1/(\alpha-1)} = \left\{
\begin{array}{ll}
x \ & \mbox{if $x=y$}, \\
\displaystyle \left( \frac{x^\alpha-y^\alpha}{\alpha(x-y)} \right)^{1/(\alpha-1)} \ & \mbox{otherwise},
\end{array}
\right.
\label{Sxy}
\end{equation}
here, $\alpha\in \mathbb{R}\backslash\{0,1\}$.
In fact, $S_\alpha(x,y)$ is known as the generalization of the logarithmic mean~\cite{L74}
\begin{equation}
\mbox{LogMean}(x,y) = \left\{
\begin{array}{ll}
x \ & \mbox{if $x=y$}, \\
\displaystyle \frac{x-y}{\ln x - \ln y} \ & \mbox{otherwise}.
\end{array}
\right.
\label{LogMean}
\end{equation} 
For given values of $\alpha$, $S_\alpha(x,y)$ reproduces known means including the logarithmic mean,
when $\alpha\to 0$, and some cases of the power mean~\cite{B03,S09}
\begin{equation}
PM_\alpha(x,y)=\left( \frac{x^\alpha+y^\alpha}{2}\right)^{1/\alpha} \, .
\label{Mxy}
\end{equation}
As examples, in Table \ref{TableSxy} we show
some expressions for $S_\alpha(x,y)$ for selected values of $\alpha$ with their corresponding
names, when available.

\begin{table}[ht]
\caption{Expressions for the Stolarsky mean $S_\alpha(x,y)$ for selected values of $\alpha$.} 
\centering 
\begin{tabular}{r l l} 
\hline\hline 
$\alpha$ & $S_\alpha(x,y)$ & name (when available) \\ [0.5ex] 
\hline 
$-\infty$ & $S_{\alpha\to-\infty}(x,y)=\min(x,y)$ & minimum value, $PM_{\alpha\to-\infty}(x,y)$  \\ [1ex]
$-4$ & $\displaystyle S_{-4}(x,y)=\left( \frac{x^3+x^2y+xy^2+y^3}{4x^4y^4} \right)^{-1/5}$ &   \\ [1ex]
$-3$ & $\displaystyle S_{-3}(x,y)=\left( \frac{x^2+xy+y^2}{3x^3y^3}\right)^{-1/4}$ &  \\ [1ex]
$-2$ & $\displaystyle S_{-2}(x,y)=\left( \frac{x+y}{2x^2y^2}\right)^{-1/3}$ &  \\ [1ex]
$-1$ & $S_{-1}(x,y)=\sqrt{xy}$ & geometric mean, $PM_{\alpha\to 0}(x,y)$  \\ [1ex]
0 & $S_{\alpha\to 0}(x,y)=\left\{
\begin{array}{ll}
x \ & \mbox{if $x=y$} \\
\displaystyle \frac{x-y}{\ln x - \ln y} \ & \mbox{otherwise}
\end{array}\right.$ & $\mbox{LogMean}(x,y)$  \\ [1ex]
$1/2$ & $\displaystyle S_{1/2}(x,y)=\left( \frac{\sqrt{x}+\sqrt{y}}{2}\right)^2$ & $PM_{1/2}(x,y)$  \\ [1ex]
1 & $S_{\alpha\to 1}(x,y)=\left\{
\begin{array}{ll}
x \ & \mbox{if $x=y$} \\
\displaystyle \frac{x-y}{x\ln x - y\ln y} \ & \mbox{otherwise}
\end{array}\right.$ & identric mean  \\ [1ex]
2 & $\displaystyle S_2(x,y)=\frac{x+y}{2}$ & arithmetic mean, $PM_1(x,y)$  \\ [1ex]
3 & $\displaystyle S_3(x,y)=\left( \frac{x^2+xy+y^2}{3}\right)^{1/2}$ &   \\ [1ex]
4 & $\displaystyle S_4(x,y)=\left( \frac{x^3+x^2y+xy^2+y^3}{4}\right)^{1/3}$ &   \\ [1ex]
$\infty$ & $S_{\alpha\to\infty}(x,y)=\max(x,y)$ & maximum value, $PM_{\alpha\to\infty}(x,y)$  \\ [1ex] 
\hline 
\end{tabular}
\label{TableSxy} 
\end{table}

Also, there is a well-known inequality relating the Stolarsky mean and the power mean, namely~\cite{L74,OT57,C66}:
\begin{equation}
S_{-1}(x,y) = PM_{\alpha\to 0}(x,y) \le S_{\alpha\to 0}(x,y) \le PM_{1/3}(x,y) \le S_2(x,y) = PM_1(x,y)
\label{ineq}
\end{equation}
or more explicitely
$$
\sqrt{xy} \le \mbox{LogMean}(x,y) \le \left( \frac{x^{1/3}+y^{1/3}}{2} \right)^3 \le \frac{x+y}{2} \, ,
$$
where the equality is attained when $x=y$.

\section{Stolarsky--Puebla index}

A large number of graph invariants of the form 
\begin{equation}
TI(G) = \sum_{uv \in E(G)} F(d_u,d_v)
\label{TI}
\end{equation}
are currently been studied in 
mathematical chemistry; where $uv$ denotes the edge of the graph $G$ connecting the vertices $u$ and $v$,
$d_u$ is the degree of the vertex $u$, and $F(x,y)$ is an appropriate chosen function, see e.g.~\cite{G13}. 

Inspired by the Stolarsky mean and given a simple graph $G=(V(G),E(G))$, here we choose 
the function $F(x,y)$ in Eq.~(\ref{TI}) as the Stolarsky mean $S_\alpha(x,y)$ and define the degree--based 
variable topological index 
\begin{equation}
SP_\alpha(G) =  S_\alpha(d_u,d_v) = \sum_{uv \in E(G)}  \left\{
\begin{array}{ll}
d_u \ & \mbox{if $d_u=d_v$}, \\
\displaystyle \left( \frac{d_u^\alpha-d_v^\alpha}{\alpha(d_u-d_v)} \right)^{1/(\alpha-1)} \ & \mbox{otherwise},
\end{array}
\right.
\label{SG}
\end{equation}
where $uv$ denotes the edge of the graph $G$ connecting the vertices $u$ and $v$,
$d_u$ is the degree of the vertex $u$, and $\alpha\in \mathbb{R}\backslash\{0,1\}$. 
We name $SP_\alpha(G)$ as the Stolarsky--Puebla index.

Note, that for given values of $\alpha$, $SP_\alpha(G)$ is related to widely studied topological indices: 
$SP_{-1}(G) = R^{-1}(G)$, where $R^{-1}(G)$ is the reciprocal Randic index~\cite{GFE14}, 
$SP_{1/2}(G) = 2^{-2} KA^1_{1/2,2}(G)$, where $KA^1_{\alpha,\beta}(G)$ is the first $(\alpha,\beta)-KA$ index~\cite{K19}, and 
$SP_2(G) = M_1(G)/2$, where $M_1(G)$ is the first Zagreb index~\cite{GT72}.
Also, for selected values of $\alpha$, $SP_\alpha(G)$ reproduces several mean Sombor indices
\begin{equation}
mSO_\alpha(G)= \sum_{uv \in E(G)} \left( \frac{d_u^\alpha+d_v^\alpha}{2}\right)^{1/\alpha} ;
\label{MG}
\end{equation}
recently introduced in~\cite{AMMS21}.
In Table \ref{TableSG} we report some expressions for $SP_\alpha(G)$ for selected values of $\alpha$ 
that we identify with known topological indices, when applicable.

\begin{table}[ht]
\caption{Expressions for the Stolarsky--Puebla index $SP_\alpha(G)$ for selected values of $\alpha$.} 
\centering 
\begin{tabular}{r l l} 
\hline\hline 
$\alpha$ & $SP_\alpha(G)$ & index equivalence \\ [0.5ex] 
\hline 
$-\infty$ & $\displaystyle SP_{\alpha\to-\infty}(G)=\sum_{uv \in E(G)} \min(d_u,d_v)$ &  $mSO_{\alpha\to-\infty}(G)$ \\ [1ex]
$-4$ & $\displaystyle SP_{-4}(G)= \sum_{uv \in E(G)} \left( \frac{d_u^3+d_u^2d_v+d_ud_v^2+d_v^3}{4d_u^4d_v^4} \right)^{-1/5}$ &   \\ [1ex]
$-3$ & $\displaystyle SP_{-3}(G)= \sum_{uv \in E(G)} \left( \frac{d_u^2+d_ud_v+d_v^2}{3d_u^3d_v^3}\right)^{-1/4}$ &  \\ [1ex]
$-2$ & $\displaystyle SP_{-2}(G)= \sum_{uv \in E(G)} \left( \frac{d_u+d_v}{2d_u^2d_v^2}\right)^{-1/3}$ &  \\ [1ex]
$-1$ & $\displaystyle SP_{-1}(G)= \sum_{uv \in E(G)} \sqrt{d_ud_v}$ & $R^{-1}(G)=mSO_{\alpha\to 0}(G)$  \\ [1ex]
0 & $\displaystyle SP_{\alpha\to 0}(G)= \sum_{uv \in E(G)} \left\{
\begin{array}{ll}
d_u \ & \mbox{if $d_u=d_v$} \\
\displaystyle \frac{d_u-d_v}{\ln d_u - \ln d_v} \ & \mbox{otherwise}
\end{array}\right.$ & logarithmic--mean index, see Eq.~(\ref{LogMeanG})  \\ [1ex]
$1/2$ & $\displaystyle SP_{1/2}(G)= \sum_{uv \in E(G)} \left( \frac{\sqrt{d_u}+\sqrt{d_v}}{2}\right)^2$ & $2^{-2} KA^1_{1/2,2}(G)=mSO_{1/2}(G)$   \\ [1ex]
1 & $\displaystyle SP_{\alpha\to 1}(G)= \sum_{uv \in E(G)} \left\{
\begin{array}{ll}
d_u \ & \mbox{if $d_u=d_v$} \\
\displaystyle \frac{d_u-d_v}{d_u\ln d_u - d_v\ln d_v} \ & \mbox{otherwise}
\end{array}\right.$ &  identric--mean index, see Eq.~(\ref{idLogMeanG}) \\ [1ex]
2 & $\displaystyle SP_2(G)= \sum_{uv \in E(G)} \frac{d_u+d_v}{2}$ & $2^{-1}M_1(G)=mSO_1(G)$  \\ [1ex]
3 & $\displaystyle SP_3(G)= \sum_{uv \in E(G)} \left( \frac{d_u^2+d_ud_v+d_v^2}{3}\right)^{1/2}$ &   \\ [1ex]
4 & $\displaystyle SP_4(G)= \sum_{uv \in E(G)} \left( \frac{d_u^3+d_u^2d_v+d_ud_v^2+d_v^3}{4}\right)^{1/3}$ &   \\ [1ex]
$\infty$ & $\displaystyle SP_{\alpha\to\infty}(G)= \sum_{uv \in E(G)} \max(d_u,d_v)$ &  $mSO_{\alpha\to\infty}(G)$ \\ [1ex] 
\hline 
\end{tabular}
\label{TableSG} 
\end{table}

\section{Computational study of $SP_\alpha(G)$ on random networks}

As a first test of the Stolarsky--Puebla index, here we apply it on two models of random networks: 
Erd\"os-R\'enyi (ER) networks and random geometric (RG) graphs. 
ER networks~\cite{SR51,ER59,ER60a,ER60b} $G_{\tbox{ER}}(n,p)$ are formed by $n$ vertices connected independently 
with probability $p \in [0,1]$. 
While RG graphs~\cite{DC02,P03} $G_{\tbox{RG}}(n,r)$ consist of $n$ vertices uniformly and independently 
distributed on the unit square, where two vertices are connected by an edge if their Euclidean distance is less 
or equal than the connection radius $r \in [0,\sqrt{2}]$.

We stress that the computational study of the Stolarsky--Puebla index we perform here is justified by the 
random nature of the network models we want to explore. Since a given parameter set [$(n,p)$ or $(n,r)$] 
represents an infinite-size ensemble of random [ER or RG] networks, the computation of $SP_\alpha(G)$ 
on a single network is irrelevant. In contrast, the computation of the average value of $SP_\alpha(G)$ on
a large ensemble of random networks, all characterized by the same parameter set, may provide useful 
{\it average} information about the full ensemble. This {\it statistical} approach, well known in random matrix 
theory studies, has been recently applied to random networks by means of topological indices, 
see e.g.~\cite{MMRS20,AHMS20,MMRS21}. Moreover, it has been shown that average topological indices
may serve as complexity measures equivalent to standard random matrix theory measures~\cite{AMRS20,AMRS21}.

\subsection{$SP_\alpha(G)$ on Erd\"os-R\'enyi random networks}
\label{ER}

In what follows we present the average values of selected Stolarsky--Puebla indices. 
All averages are computed over ensembles of $10^7/n$ ER networks characterized 
by the parameter pair $(n,p)$.

\begin{figure}[t!]
\begin{center} 
\includegraphics[width=0.65\textwidth]{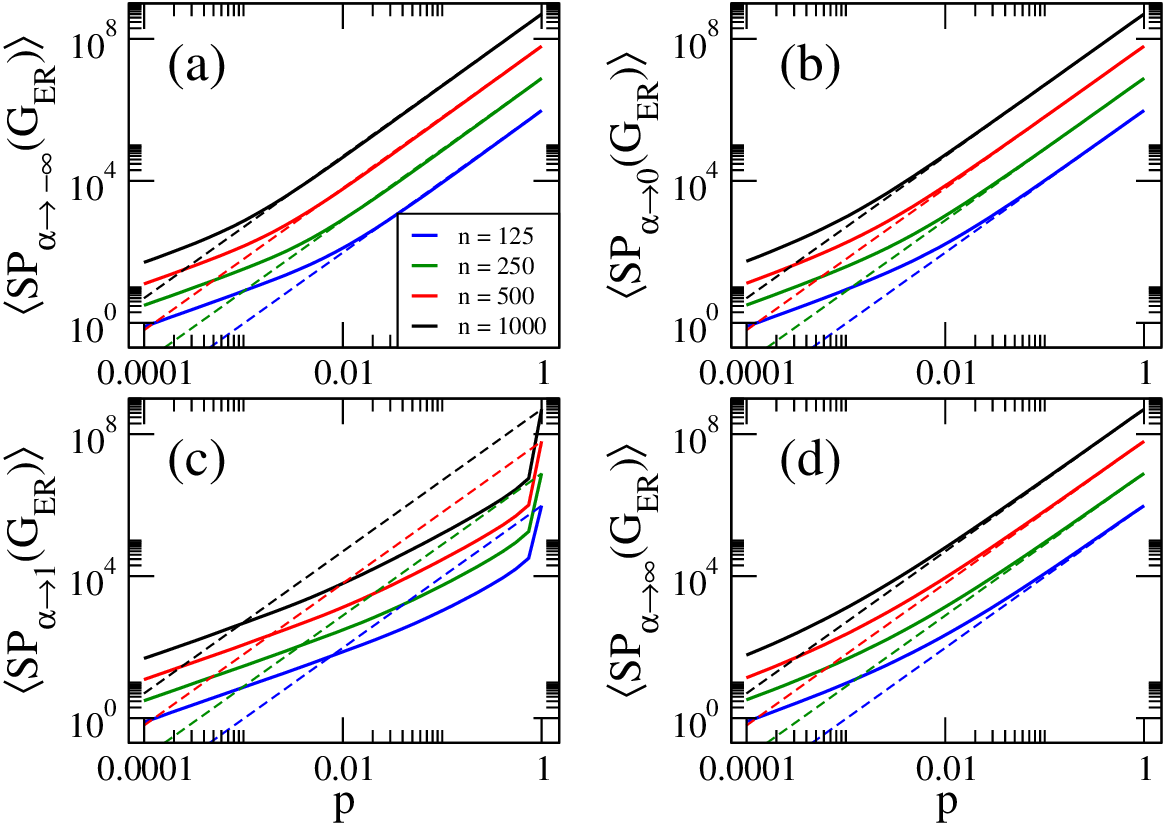}
\caption{\footnotesize{
Average value of the Stolarsky--Puebla index $\left< SP_\alpha(G_{\tbox{ER}}) \right>$ as a function of the 
probability $p$ of Erd\"os-R\'enyi networks of size $n$. Here
(a) $\alpha\to -\infty$,
(b) $\alpha\to 0$,
(c) $\alpha\to 1$, and
(d) $\alpha\to \infty$.
Dashed lines correspond to Eq.~(\ref{SofER}).}}
\label{Fig01}
\end{center}
\end{figure}

In Fig.~\ref{Fig01} we present the average Stolarsky--Puebla index $\left< SP_\alpha(G_{\tbox{ER}}) \right>$
for $\alpha\to -\infty$, $\alpha\to 0$, $\alpha\to 1$, and $\alpha\to \infty$ as a function of the 
probability $p$ of ER networks of sizes 
$n=\{125,250,500,1000\}$. From this figure we observe that the curves of 
$\left< SP_\alpha(G_{\tbox{ER}}) \right>$ are monotonically increasing functions of $p$.

We note that in the dense limit, i.e. when $np\gg 1$, we can approximate 
$d_u \approx d_v \approx \left<  d \right>$ in Eq.~(\ref{SG}), with
\begin{equation}
\label{k}
\left< d \right> = (n-1)p .
\end{equation}
Thus, when $np\gg 1$, we can approximate $SP_\alpha(G_{\tbox{ER}})$ as
\begin{equation}
\label{SofER}
SP_\alpha(G_{\tbox{ER}}) \approx \sum_{uv \in E(G)} d_u \approx \sum_{uv \in E(G)} \left<  d \right> 
\approx \frac{1}{2} n \left[  (n-1)p \right]^2 ,
\end{equation}
where we have used $|E(G_{\tbox{ER}})|=n(n-1)p/2$.
In Fig.~\ref{Fig01}, we show that Eq.~(\ref{SofER}) (dashed lines) indeed describes well the data 
(thick full curves) for large enough $p$; except for the case $\left< SP_{\alpha\to 1}(G_{\tbox{ER}}) \right>$,
see Fig.~\ref{Fig01}(c).
We also verified that Eq.~(\ref{SofER}) describes well the data for other values of $\alpha$,
however we did not include them in Fig.~\ref{Fig01} to avoid figure saturation.
We also observed that the smaller the value of $\alpha$ the wider the range of $p$ where the 
coincidence between Eq.~(\ref{SofER}) and the computational data is observed; compare for example 
Figs.~\ref{Fig01}(a) and~\ref{Fig01}(d), where it is clear that the correspondence of the computational
data with Eq.~(\ref{SofER}) is much better in the case of $\alpha\to -\infty$ than for $\alpha\to \infty$.
In addition, it is relevant to note that Eq.~(\ref{SofER}) does not depend on $\alpha$.

We also notice that in Fig.~\ref{Fig01} we present average Stolarsky--Puebla indices as a function of the 
probability $p$ of ER networks of four different sizes $n$. 
It is quite clear from these figures that the curves, characterized by the different network sizes, 
are very similar but displaced on both axes. 
This behavior suggests that the average Stolarsky--Puebla indices can be scaled, as will be shown below.

\begin{figure}[h!]
\begin{center} 
\includegraphics[width=0.65\textwidth]{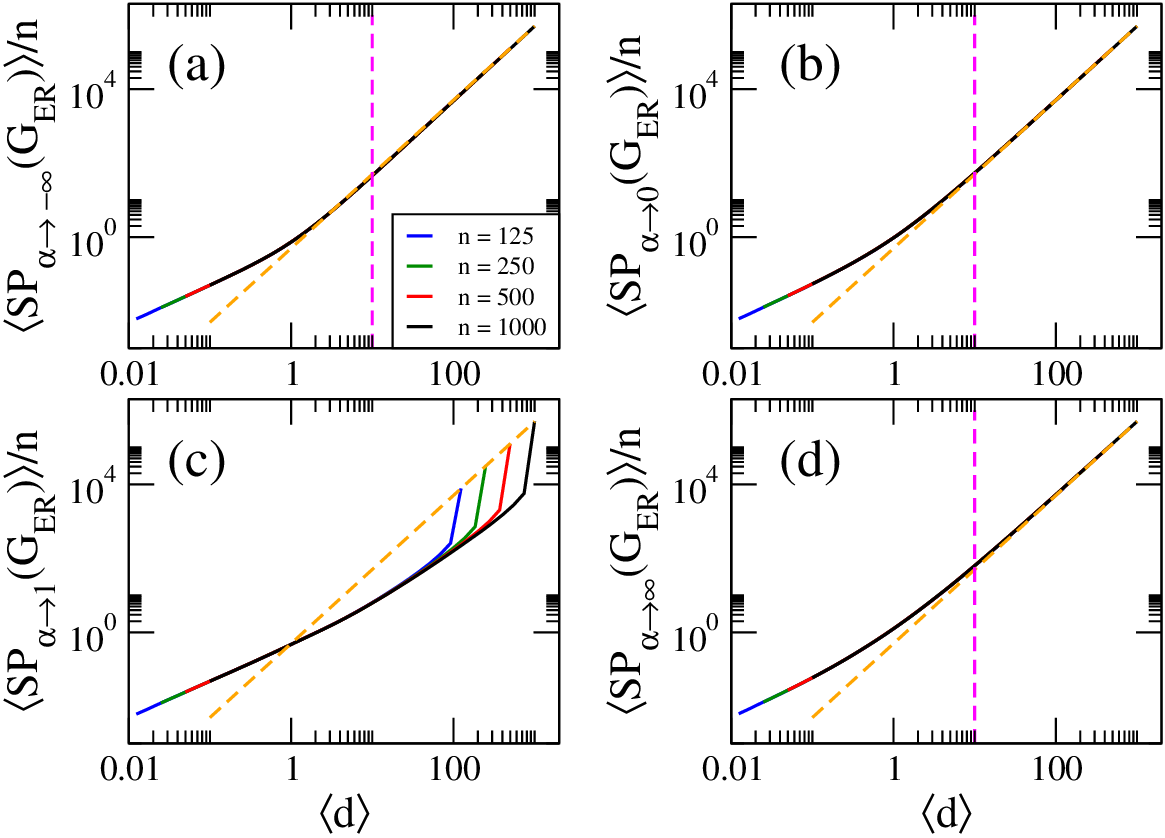}
\caption{\footnotesize{
Average value of the Stolarsky--Puebla index $\left< SP_\alpha(G_{\tbox{ER}}) \right>$, normalized to the network 
size $n$, as a function of the average degree $\left< d \right>$ of Erd\"os-R\'enyi networks. 
Same curves as in Fig.~\ref{Fig01}.
Orange dashed lines are Eq.~(\ref{SnofER}).
The vertical magenta dashed lines indicate $\left< d \right>=10$.}}
\label{Fig02}
\end{center}
\end{figure}

From Eq.~(\ref{SofER}) we observe that $\left< SP_\alpha(G_{\tbox{ER}}) \right>\propto n f[(n-1)p)]$ or
\begin{equation}
\label{scaling}
\left< SP_\alpha(G_{\tbox{ER}}) \right>\propto n f(\left< d \right>).
\end{equation}
Therefore, in Fig.~\ref{Fig02} we plot again the average Stolarsky--Puebla indices reported in Fig.~\ref{Fig01},
but now normalized to $n$, as a function of $\left< d \right>$ showing that all indices are now properly 
scaled; i.e.~the curves painted in different colors for different network sizes fall on top of each other. 
Moreover, we can rewrite Eq.~(\ref{scaling}) as
\begin{equation}
\label{SnofER}
\frac{\left< SP_\alpha(G_{\tbox{ER}}) \right>}{n} \approx \frac{1}{2} \left<  d \right>^2 .
\end{equation}
In Fig.~\ref{Fig02}, we show that Eq.~(\ref{SnofER}) (orange-dashed 
lines) indeed describe well the data (thick full curves) for $\left< d \right>\ge 10$;
except for $\left< SP_{\alpha\to 1}(G_{\tbox{ER}}) \right>$, see Fig.~\ref{Fig02}(c).

It is relevant to stress that even when Eq.~(\ref{scaling}) was expected to be valid in the dense limit 
(i.e.~for $\left< d \right> \gg 1$), it is indeed valid for any $\left< d \right>$
as clearly seen in Fig.~\ref{Fig02}.

\subsection{$SP_\alpha(G)$ on random geometric graphs}
\label{RG}

As in the previous Subsection, here we present the average values of selected Stolarsky--Puebla indices.
Again, all averages are computed over ensembles of $10^7/n$ 
random graphs, each ensemble characterized by a fixed parameter pair $(n,r)$. 

\begin{figure}[t!]
\begin{center} 
\includegraphics[width=0.65\textwidth]{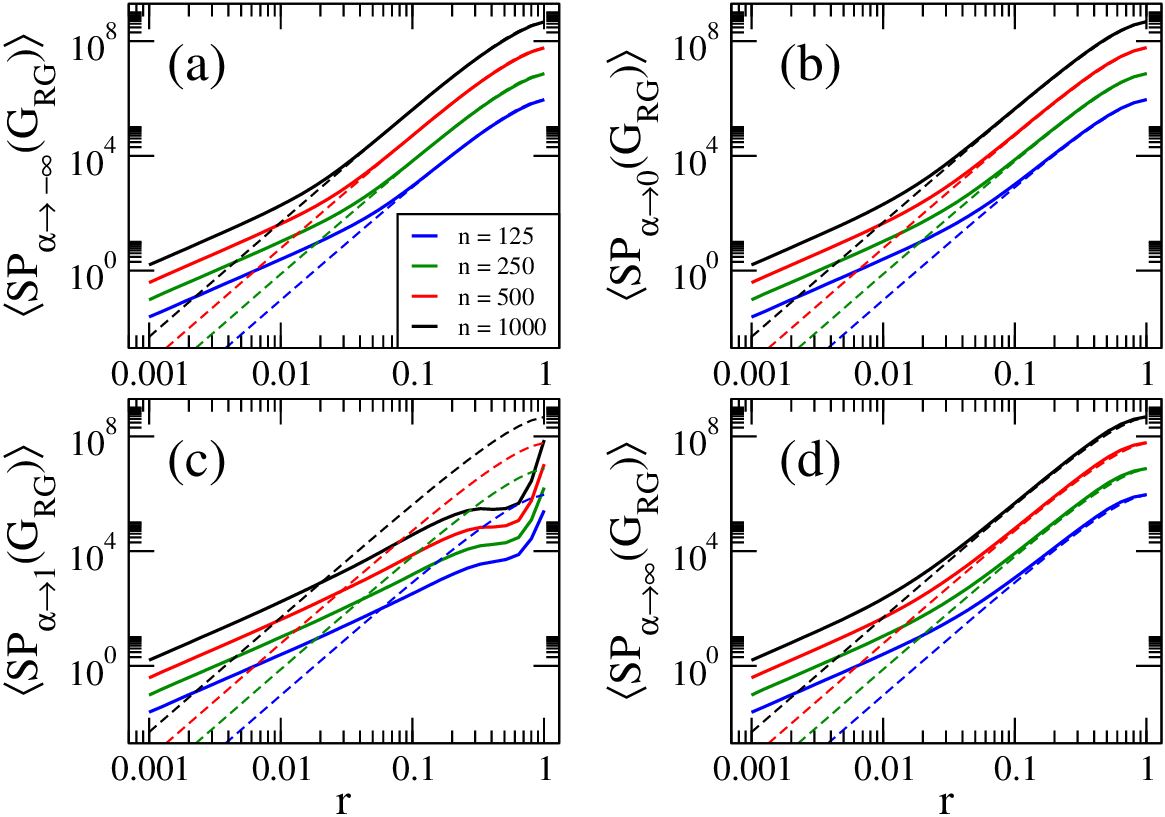}
\caption{\footnotesize{
Average value of the Stolarsky--Puebla index $\left< SP_\alpha(G_{\tbox{RG}}) \right>$ as a function 
of the connection radius $r$ of random geometric graphs of size $n$. Here
(a) $\alpha\to -\infty$,
(b) $\alpha\to 0$,
(c) $\alpha\to 1$, and
(d) $\alpha\to \infty$.
Dashed lines correspond to Eq.~(\ref{SofRG}).}}
\label{Fig03}
\end{center}
\end{figure}

In Fig.~\ref{Fig03} we present the average Stolarsky--Puebla index $\left< SP_\alpha(G_{\tbox{ER}}) \right>$
for $\alpha\to -\infty$, $\alpha\to 0$, $\alpha\to 1$, and $\alpha\to \infty$ as a function of the connection 
radius $r$ of RG graphs of sizes $n=\{125,250,500,1000\}$.
For comparison purposes, Fig.~\ref{Fig03} is equivalent to Fig.~\ref{Fig01}. In fact, all the observations 
made in the previous Subsection for ER networks are also valid for RG graphs by just replacing 
$G_{\tbox{ER}}\to G_{\tbox{RG}}$ and $p\to g(r)$, with~\cite{EM15}
\begin{equation}
g(r) = 
\left\{ 
\begin{array}{ll}
           r^2  \left[ \pi - \frac{8}{3}r +\frac{1}{2}r^2 \right] & \quad 0 \leq r \leq 1 \, , 
           \vspace{0.25cm} \\
             \frac{1}{3} - 2r^2 \left[ 1 - \arcsin(1/r) + \arccos(1/r) \right] 
               +\frac{4}{3}(2r^2+1) \sqrt{r^2-1} 
             -\frac{1}{2}r^4 & \quad 1 \leq r \leq \sqrt{2} \, .
\end{array}
\right.
\label{g(r)}
\end{equation}

As well as for ER networks, here, in the dense limit, when $nr\gg 1$, we can approximate 
$d_u \approx d_v \approx \left<  d \right>$ with
\begin{equation}
\label{kRG}
\left< d \right> = (n-1)g(r) .
\end{equation}
Therefore, in the dense limit, $SP_\alpha(G_{\tbox{RG}})$ is well approximated by:
\begin{equation}
\label{SofRG}
SP_\alpha(G_{\tbox{RG}}) \approx  \frac{1}{2} n \left[  (n-1)g(r) \right]^2 .
\end{equation}
In Fig.~\ref{Fig03}, we show that Eq.~(\ref{SofRG}) (dashed 
lines) indeed describes well the data (thick full curves) for large enough $r$; 
except for the case $\left< SP_{\alpha\to 1}(G_{\tbox{RG}}) \right>$,
see Fig.~\ref{Fig03}(c).

\begin{figure}[h!]
\begin{center} 
\includegraphics[width=0.65\textwidth]{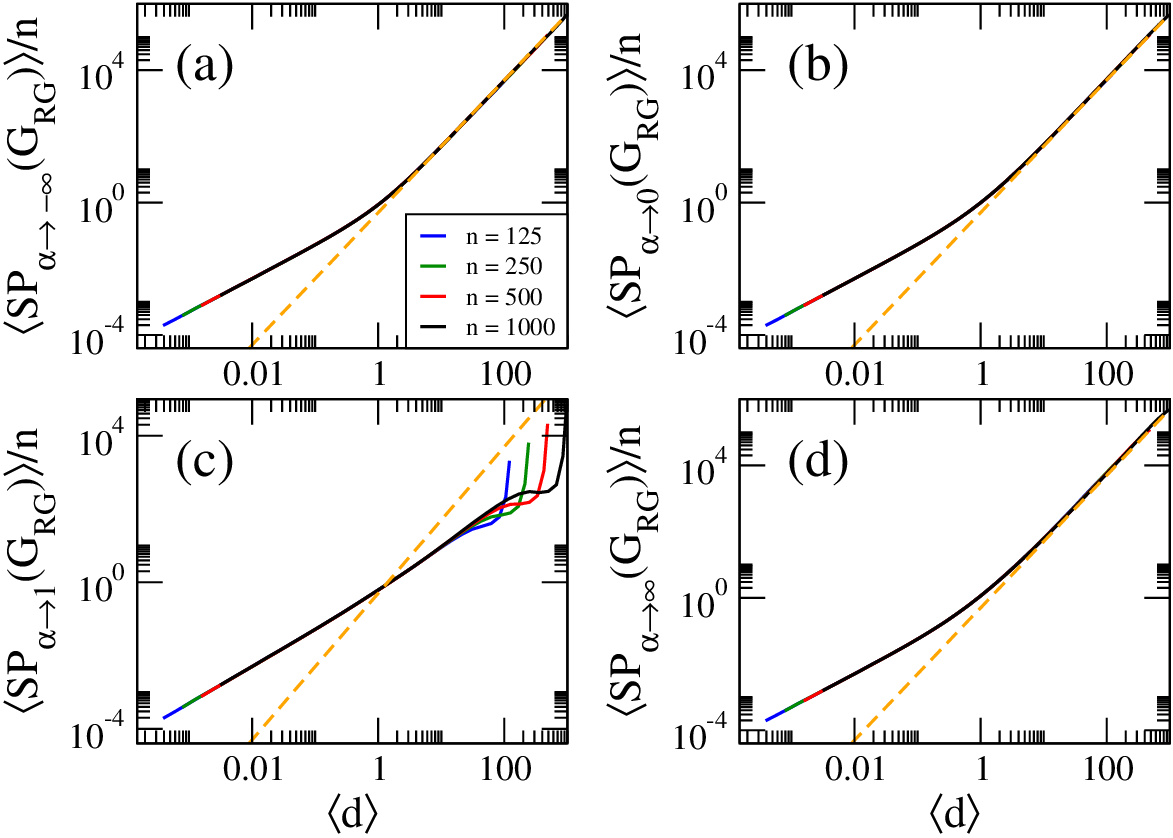}
\caption{\footnotesize{
Average value of the Stolarsky--Puebla index $\left< SP_\alpha(G_{\tbox{GR}}) \right>$, normalized to the network 
size $n$, as a function of the average degree $\left< d \right>$ of random geometric graphs. 
Same curves as in Fig.~\ref{Fig03}.
Orange dashed lines are Eq.~(\ref{SnofRG}).
The vertical magenta dashed lines indicate $\left< d \right>=10$.}}
\label{Fig04}
\end{center}
\end{figure}

It is quite remarkable to note that by substituting the average degree of Eq.~(\ref{kRG}) into 
Eq.~(\ref{SofRG}) we get exactly the same expression of Eq.~(\ref{SnofER}):
\begin{equation}
\label{SnofRG}
\frac{\left< SP_\alpha(G_{\tbox{RG}}) \right>}{n} \approx \frac{1}{2} \left<  d \right>^2 .
\end{equation}
So, in Fig.~\ref{Fig04} we plot again the average Stolarsky--Puebla indices reported in Fig.~\ref{Fig03} for 
RG graphs, but now normalized to $n$, as a function of $\left< d \right>$ showing that all curves 
are now properly scaled.
Also, in Fig.~\ref{Fig04}, we show that Eq.~(\ref{SnofRG}) (orange-dashed 
lines) indeed describes well the data (thick full curves) for $\left< d \right>\ge 10$.
We note that as well as for ER networks, here for RG graphs we do not observe the scaling of
$\left< SP_{\alpha\to 1}(G_{\tbox{RG}}) \right>$.

\section{Discussion and conclusions}

We have introduced a degree--based variable topological index inspired on the Stolarsky mean,
known as the generalization of the logarithmic mean.
We named this new index as the Stolarsky--Puebla index $SP_\alpha(G)$, see Eq.~(\ref{SG}).
For given values of $\alpha$, the Stolarsky--Puebla index is related to well-known topological indices, in 
particular it reproduces several mean Sombor indices $mSO_\alpha(G)$, see Eq.~(\ref{MG}).

We want to add that the inequality of Eq.~(\ref{ineq}) can be straightforwardly used to
state inequalities for the indices $SP_\alpha(G)$ and $mSO_\alpha(G)$, as well as for related indices: 
\begin{equation}
SP_{-1}(G) = mSO_{\alpha\to 0}(G) \le SP_{\alpha\to 0}(G) \le mSO_{1/3}(G) \le SP_2(G) = mSO_1(G)
\label{ineq21}
\end{equation}
or
\begin{equation}
R^{-1}(G) \le \mbox{LogMean}(G) \le mSO_{1/3}(G) \le 2^{-1}M_1(G) \, ,
\label{ineq22}
\end{equation}
which sets bounds for the logarithmic--mean topological index
\begin{equation}
\mbox{LogMean}(G) = \left\{
\begin{array}{ll}
d_u \ & \mbox{if $d_u=d_v$}, \\
\displaystyle \frac{d_u-d_v}{\ln d_u - \ln d_v} \ & \mbox{otherwise},
\end{array}
\right.
\label{LogMeanG}
\end{equation}
with respect to the reciprocal Randic index, the mean Sombor index with $\alpha=1/3$, and the first 
Zagreb index.
 
Since there are not many degree--based topological indices including logarithmic functions (as 
well-known exceptions we can mention the logarithms of the three multiplicative Zagreb
indices~\cite{G13} and the Adriatic indices~\cite{V10,V11}) we want to highlight the release of the 
logarithmic--mean topological index $\mbox{LogMean}(G)$ of Eq.~(\ref{LogMeanG}) as well as the identric--mean index 
\begin{equation}
\mbox{idLogMean}(G) = \left\{
\begin{array}{ll}
d_u \ & \mbox{if $d_u=d_v$}, \\
\displaystyle \frac{d_u-d_v}{d_u\ln d_u - d_v\ln d_v} \ & \mbox{otherwise},
\end{array}
\right.
\label{idLogMeanG}
\end{equation}
corresponding to $SP_{\alpha\to 0}(G)$ and $SP_{\alpha\to 1}(G)$, respectively.

We have also applied the Stolarsky--Puebla index $SP_\alpha(G)$ to Erd\"os-R\'enyi (ER) networks and random 
geometric (RG) graphs and within a statistical random matrix theory approach we demonstrated that 
$\left< SP_\alpha(G) \right>$, normalized to the order of the network, scales with 
the corresponding average degree $\left< d \right>$. However, it is fair to recognize that, for both random
network models, $\left< SP_{\alpha\to 1}(G) \right> = \left< \mbox{idLogMean}(G) \right>$ did not scale; 
so we believe that the identric--mean index deserves further investigation.

In addition, from Eq.~(\ref{ineq21}) we are able to write an equivalent inequality but for the corresponding 
average values:
\begin{equation}
\left< SP_{-1}(G) \right> \le \left< \mbox{LogMean}(G) \right> \le \left< mSO_{1/3}(G) \right> \le \left< SP_2(G) \right> .
\label{ineq21av}
\end{equation}
Indeed, we verified that (\ref{ineq21av}) is satisfied for both, ER random networks and RG graphs (not
shown here). Moreover, we computationally found that
\begin{equation}
\left< \mbox{idLogMean}(G) \right> \le \left< SP_{\alpha\ne 1}(G) \right> ,
\label{ineq3}
\end{equation}
for the two random network models we study here (not explicitelly shown here but partially observed in 
Figs.~\ref{Fig01} and~\ref{Fig03}).
The equalities in Eqs.~(\ref{ineq21av}) and~(\ref{ineq3}) are attained when $p=1$ and $r=\sqrt{2}$, for 
ER random networks and RG graphs, respectively.

Finally, we want to recall that through a quantitative structure property relationship (QSPR) analysis it was
shown~\cite{AMMS21} that $mSO_{\alpha\to\pm\infty}(G)$ are good predictors of the standard enthalpy of 
vaporization, the enthalpy of vaporization, and the heat of vaporization at 25$^{\circ}$C of octane isomers. 
Furthermore, since $SP_{\alpha\to\pm\infty}(G)=mSO_{\alpha\to\pm\infty}(G)$, we can conclude that 
$SP_{\alpha\to\pm\infty}(G)$ correlate well with the aforementioned physicochemical properties of octane isomers.

In future works we plan to explore mathematical and computational properties of $SP_\alpha(G)$, 
as well as finding optimal bounds and new relationships with known topological indices.

\section*{Acknowledgment}

J.A.M.-B. acknowledges financial support from CONACyT (Grant No.~A1-S-22706) and
BUAP (Grant No.~100405811-VIEP2021).
The research of J.M.S. was supported by a grant from Agencia Estatal de Investigaci\'on 
(PID2019-106433GBI00
/AEI/10.13039/501100011033), Spain.


\end{document}